\documentclass{article}

\usepackage{amsthm}
\usepackage{amsmath}
\usepackage{amssymb}
\usepackage{graphicx}
\usepackage{epsfig}

\newtheorem{theorem}{Theorem}

\newtheorem{definition}[theorem]{Definition}
\newtheorem{example}[theorem]{Example}
\newtheorem{remark}[theorem]{Remark}
\newtheorem{corollary}[theorem]{Corollary}

\begin{document}

\begin{center}
\LARGE{\textbf{INDUCED CONNECTIONS ON SUBMANIFOLDS IN SPACES WITH
FUNDAMENTAL GROUPS}}

\vspace*{3mm}

{\large  Maks A. Akivis, Vladislav V. Goldberg, and Arto V.
Chakmazyan}

\end{center}

Abstract. {\footnotesize The authors establish a relation of the
theory of varieties with degenerate Gauss maps in projective
spaces with the theory of congruences and pseudocongruences of
subspaces and show how these two theories can be applied to the
construction of induced connections on submanifolds of projective
spaces and other spaces endowed with a projective structure.}

\vspace*{2mm}


\vspace*{2mm}


 \vspace*{3mm}

\setcounter{equation}{0}

\setcounter{section}{-1}

\section{Introduction}

The theory of congruences and pseudocongruences of subspaces of a
projective space is closely related to the theory of varieties
with degenerate Gauss maps.

In a three-dimensional space $\mathbb{P}^3$ as well as in
three-dimensional spaces endowed with a projective structure (such
as an affine, Euclidean, and non-Euclidean space), the theory of
congruences was studied by many geometers. The extensive
monographs on this subject were published (see, for example, the
monograph [Fi 50]).

In this paper, we establish a relation of the theory of varieties
with degenerate Gauss maps in projective spaces with the theory of
congruences and pseudocongruences of subspaces and show how these
two theories can be applied to the construction of induced
connections on submanifolds of projective spaces and other spaces
endowed with a projective structure.


\section{Basic Equations of a Variety with  \\ a Degenerate Gauss Map.}

A smooth $n$-dimensional variety $X$ of a projective space
$\mathbb{P}^N$ is called  {\em tangentially degenerate} or {\em a
variety with a degenerate Gauss map} if the rank of its Gauss map
$$
\gamma: X \rightarrow \mathbb{G} (n, N)
$$
is less than $n, \; 0 \leq r = \mbox{{\rm rank}} \; \gamma < n$.
Here $x \in X, \; \gamma (x) = T_x (X)$, and $T_x (X)$ is the
tangent subspace to $X$ at $x$ considered as an $n$-dimensional
projective space $\mathbb{P}^n$.  The number  $r$ is also called
the {\em rank} of $X, \; r = \mbox{{\rm rank}} \; X$. The case $r
= 0 $ is trivial one: it gives just an $n$-plane.

Let $X \subset  \mathbb{P}^N$ be  an almost everywhere smooth
$n$-dimensional variety with a degenerate Gauss map. Suppose that
$0 < \mbox{{\rm rank}} \, \gamma = r < n$. Denote by $L$ a leaf of
this map, $L = \gamma^{-1} (T_x) \subset X;\, \dim L = n - r = l$.

As was proved in [AG 04] (see Theorem 3.1 on p. 95), a variety
with a degenerate Gauss map of rank $r$ foliates into its plane
leaves $L$ of dimension $l$, along which the tangent subspace $T_x
(X)$ is fixed.

The tangent subspace $T_x (X)$ is fixed when a point $x$ moves
along  regular points of $L$. This is the reason that we denote it
by $T_L, \, L \subset T_L$. A pair $(L, T_L)$ on $X$ depends on
$r$ parameters.

The foliation on $X$ defined as indicated above is called the {\em
Monge--Amp\`{e}re foliation}.

The varieties of rank $r < n$ are multidimensional analogues of
developable surfaces of a three-dimensional Euclidean space.

The main  results on the geometry of varieties with degenerate
Gauss maps and further references can be found in Chapter 4 of the
book [AG 93] and in the recently published book [AG 04].

In this section, we find the basic equations of a variety $X$ with
a degenerate Gauss map of dimension $n$ and rank $r$ in a
projective space $\mathbb{P}^N$.

In what follows, we will use the following ranges of indices:
$$
a, b, c = 1, \ldots, l; \; p, q = l + 1, \ldots , n; \; \alpha,
\beta = n + 1, \ldots , N.
$$

A point $x \in X$ is said to be a \textit{regular point} of the
map $\gamma$ and of the variety $X$ if $\dim T_x X = \dim X = n$,
and a point $x \in X$ is called \textit{singular} if $\dim T_x X
> \dim X = n$.

We associate a family of moving frames $\{A_u\}, \, u = 0, 1,
\ldots , N$, with $X$ in such a way that the point $A_0 = x$ is a
regular point of $X$; the points $A_a$ belong to the leaf $L$ of
the Monge--Amp\`{e}re foliation passing through the point $A_0$;
the points $A_p$ together with the points $A_0, A_a$ define the
tangent subspace $T_L X$ to $X$; and the points $A_\alpha$ are
located outside of the subspace $T_L X$.

The equations of infinitesimal displacement of the moving frame
$\{A_u\}$ are
 \begin{equation}\label{eq:1}
 dA_u = \omega_u^v A_v, \;\;\;\;u, v = 0, 1, \ldots , N,
\end{equation}
where $\omega_u^v$ are 1-forms satisfying the structure equations
of the projective space $\mathbb{P}^N$:
 \begin{equation}\label{eq:2}
 d \omega_u^v = \omega_u^w \wedge \omega_w^v, \;\;\;\;u, v, w = 0, 1, \ldots , N.
\end{equation}

As a result of the specialization of the moving frame mentioned
above, we obtain the following equations of the variety $X$ (see
[AG 04], Section 3.1):
\begin{equation}\label{eq:3}
\omega_0^\alpha = 0,
\end{equation}
\begin{equation}\label{eq:4}
\omega_a^\alpha = 0,
\end{equation}
\begin{equation}\label{eq:5}
\omega_p^\alpha = b_{pq}^\alpha \omega^q, \;\; b_{pq}^\alpha =
b_{qp}^\alpha,
\end{equation}
\begin{equation}\label{eq:6}
\omega^p_a = c^p_{aq} \omega^q.
\end{equation}

The 1-forms $\omega^q := \omega_0^q$ in these equations are basis
forms of the Gauss image $\gamma (X)$ of the variety $X$, and the
quantities $b_{pq}^\alpha$ form the second fundamental tensor of
the variety $X$ at the point $x = A_0$. The quantities
$b_{pq}^\alpha$ and $c^p_{aq}$ are related by the following
equations:
\begin{equation}\label{eq:7}
b_{sq}^\alpha c^s_{ap} = b_{sp}^\alpha c^s_{aq}.
\end{equation}

Equations (5) and (6) are called the \emph{basic equations} of a
variety $X$ with a degenerate Gauss map (see [AG 04], Section
3.1).

Note that  under transformations of the points $A_p$, the
quantities $c^p_{aq}$ are transformed as tensors with respect to
the indices $p$ and $q$. As to the index $a$, the  quantities
$c^p_{aq}$ do not form a tensor with respect to this index.
Nevertheless, under transformations of the points $A_0$ and $A_a$,
the quantities $c^p_{aq}$ along with the unit tensor $\delta_q^p$
are transformed as  tensors. For this reason, the system of
quantities $c^p_{aq}$ is called a {\em quasitensor}.

Denote by $B^\alpha$ and $C_a$  the $(r \times r)$-matrices of
coefficients occurring in equations (5) and (6):
$$
B^\alpha = (b^\alpha_{pq}), \;\; C_a = ( c^p_{aq}).
$$
Sometimes we will use the identity matrix $C_0 = (\delta^p_q)$ and
the index \linebreak $i = 0, 1, \ldots, l$, i.e., $\{i\} = \{0,
a\}$. Then equations (5) and (7) can be combined and written as
follows:
$$
(B^\alpha C_i)^T = (B^\alpha C_i),
$$
i.e., the matrices
$$
H^\alpha_i = B^\alpha C_i = (b^\alpha_{qs} c^s_{ip})
$$
are symmetric.

The quadratic forms
\begin{equation}\label{eq:8}
\Phi^\alpha = b_{pq}^\alpha \omega^p \omega^q
\end{equation}
are the \textit{second fundamental forms} of the variety $X$ at
the point $x = A_0$, and the forms
$$
\Phi^\alpha = b^\alpha_{ps} (\delta^s_q x^0 + c^s_{aq} x^a) \,
\omega^p \, \omega^q
$$
are the \textit{second fundamental forms} of the variety $X$ at
the point $x = x^0 A_0 + x^a A_a \in L$.

Let $\{\alpha^u\}$ be the dual coframe (or tangential frame) in
the space $(\mathbb{P}^{N})^*$  to the frame $\{A_u\}$. Then the
hyperplanes $\alpha^u$ of the frame $\{\alpha^u\}$  are connected
with the points of the frame $\{ A_u\}$ by the conditions
\begin{equation}\label{eq:9}
(\alpha^u, A_v) = \delta^u_v.
\end{equation}

Conditions (9) mean that the hyperplane $\alpha^u$ contains all
points $A_v, \; v \neq u$, and that the condition of normalization
$(\alpha^u, A_u) = 1$ holds.

The equations
$$
\xi_\beta \alpha^\beta = 0
$$
defines the system of tangent hyperplanes passing through the
tangent subspace $T_L X$, and this system of tangent hyperplanes
defines the system of second fundamental forms
\begin{equation}\label{eq:10}
 II = \xi_\beta b_{pq}^\beta \omega^p \omega^q
\end{equation}
and the system of second fundamental tensors
$$
 \xi_\beta b_{pq}^\beta
$$
of the variety $X$  at the point $x = A_0$.

A variety $X$ with a degenerate Gauss map is \textit{dually
nondegenerate} if the dimension of its dual variety $X^*$ equals
$N - l -1$. By the generalized Griffiths--Harris theorem (see
Theorem 3.2 and Corollary 3.3 in [AG 04], pp. 97--99), \textit{a
variety $X$  with a degenerate Gauss map is dually nondegenerate
if and only if at any smooth point $x \in X$ there is at least one
nondegenerate second fundamental form in the system of second
fundamental forms $\xi_\alpha b_{pq}^\alpha \omega^p \omega^q$ of
$X$.}

In what follows, we will consider only dually nondegenerate
varieties $X$ with degenerate Gauss maps.

\section{Focal Images of a Variety with a Degenerate Gauss Map.}

Suppose that  $X$ is a variety with a degenerate Gauss map of
dimension $n$ and rank $r$ in the space $\mathbb{C} \mathbb{P}^N$.
As we noted earlier, such a variety carries an \linebreak
$r$-parameter family of $l$-dimensional plane generators $L$ of
dimension $l = n - r$. Let $x = x^0 A_0 + x^a A_a$ be an arbitrary
point of a generator $L$. For such a point, we have
$$
d x = (dx^0 + x^0 \omega_0^0 + x^a \omega_a^0) A_0 + (dx^a + x^0
\omega^a + x^b \omega_b^a) A_a + (x^0 \omega^p + x^a \omega_a^p)
A_p.
$$
By (6), it follows that
\begin{equation}\label{eq:11}
d x \equiv (x^0 \delta_q^p + x^a c^p_{aq}) A_p \omega^q \pmod{L}.
\end{equation}
The matrix $(J^p_q) = (x^0 \delta_q^p + x^a c^p_{aq})$ is the
\emph{Jacobi matrix} of the map $\gamma: X \rightarrow \mathbb{G}
(n, N)$, and the determinant
$$
J (x) = \det \,(J^p_q ) = \det \,(x^0 \delta_q^p + x^a c^p_{aq})
$$
of this matrix is the \emph{Jacobian} of the map $\gamma$.

It is easy to see that $J (x) \neq 0$  at regular points and  $J
(x) = 0$ at singular points.

By (6) and (11), the set of singular points of a generator $L$ of
the variety $X$ is determined by the equation
\begin{equation}\label{eq:12}
 \det \,(\delta_q^p x^0 + c^p_{aq}  x^a) = 0.
\end{equation}
Hence this set is an \emph{algebraic hypersurface of dimension $l
- 1$ and degree $r$ in the generator $L$}. This hypersurface (in
$L$) is called the {\em focus hypersurface} and  is denoted $F_L$.

Because for $x^a = 0$ the left-hand side of equation (12) takes
the form
$$
\det \,(x^0 \delta^p_q) = (x^0)^r,
$$
it follows that the point $A_0$ is a regular point of the
generator $L$.

We call a tangent hyperplane
$\xi = (\xi_\alpha)$ {\em singular} (or a \emph{focus hyperplane})
if
\begin{equation}\label{eq:13}
\det \,(\xi_\alpha b^\alpha_{pq}) = 0,
\end{equation}
i.e., if the rank of matrix $(\xi_\alpha b^\alpha_{pq})$ is
reduced. Condition (13) is an equation of degree $r$ with respect
to the tangential coordinates $\xi_\alpha$ of the hyperplanes
$\xi$ containing the tangent subspace $T_L (X)$.

Because we consider only dually nondegenerate  varieties with
degenerate Gauss maps, \textit{there exists at least one
nondegenerate form in the system of second fundamental forms of
$X$} (see the end of Section 1). Hence in the space
$\mathbb{P}^N$, equation (13) defines an algebraic hypercone of
degree $r$, whose vertex is the tangent subspace $T_L (X)$. This
hypercone is called the {\em focus hypercone} and is denoted
$\Phi_L$ (see [AG 93], p. 119).


Note that if a variety $X$ is dually degenerate, then on such a
variety, equations (13) are satisfied identically, and the variety
$X$ does not have focus hypercones.

The focus hypersurface $F_L \subset L$  and the focus hypercone
$\Phi_L$ with vertex  $T_L$ are called the {\em focal images} of
the variety $X$ with a degenerate Gauss map.

\section{Congruences and Pseudocongruences in \\ a
Projective Space}

In a projective space $\mathbb{P}^n$, we consider a family $Y$ of
its $l$-dimensional subspaces $L, \, \dim L = l$, which depends on
$r = n - l$ parameters. We assume that not more than a finite
number of subspaces $L$ passes through every point $x \in
\mathbb{P}^n$. If we restrict ourselves by a small neighborhood of
a subspace $L$, then we can assume that only one subspace $L
\subset Y$ passes through a generic point $x \in L$.  Such
families of the space $\mathbb{P}^n$ are called the
\textit{congruences}.

The dual image for a congruence $Y$ of $l$-dimensional subspaces
in $\mathbb{P}^n$ is a \textit{pseudocongruence} $Y^*$ which is an
$r$-parameter family of subspaces of dimension $r - 1$. Every
hyperplane $\xi \subset \mathbb{P}^n$ contains not more than a
finite number of subspaces $L^* \subset Y^*$. However, if we
consider an infinitesimally small neighborhood of the subspace
$L^*$ of the pseudocongruence $Y^*$, then there is only a single
subspace $L^*$ in the hyperplane $\xi$.

In what follows, we shall establish a relation of the theory of
varieties with degenerate Gauss maps in projective spaces with the
theory of congruences and pseudocongruences of subspaces and show
how these two theories can be applied to the construction of
induced connections on submanifolds of projective spaces and other
spaces endowed with a projective structure.

\vspace*{2mm}

So, consider in $\mathbb{P}^n$ a congruence $Y$ of $l$-dimensional
subspaces $L$. We associate with its element $L$ a family of
projective frames $\{A_0, A_1, \ldots , A_n\}$ chosen in such a
way that the points $A_0, A_1, \ldots , A_l$ are located in $L$,
and the points $A_{l+1}, \ldots , A_n$ are located outside of $L$.
The equations of infinitesimal displacement of such frames have
the form
\begin{equation}\label{14}
\renewcommand{\arraystretch}{1.3}
 \left\{
\begin{array}{ll}
dA_i = \omega_i^j A_j + \omega_i^p A_p, \\
dA_p = \omega_p^i A_i + \omega_p^q A_q, \\
\end{array} \right.
\renewcommand{\arraystretch}{1}
\end{equation}
where $i, j = 0, 1, \ldots, l; \, p, q = l + 1, \ldots , n$, and
$L = A_0 \wedge A_1 \ldots \wedge A_l$ is a generator of the
congruence $Y$ in question. Because this generator depends on $r$
parameters and is fixed, when $\omega_i^p = 0$, the forms
$\omega_i^p$ are expressed linearly in terms of the differentials
of these $r$ parameters or in terms of linearly independent
1-forms $\theta^p$---linear combinations of these differentials:
\begin{equation}\label{15}
  \omega_i^p = c_{iq}^p \theta^q.
\end{equation}
Under admissible linear transformations of the basis forms
$\theta^p$, the matrices $C_i = (c_{iq}^p)$ are transformed
according to the tensor law with respect to the indices $p$ and
$q$.

A point $F \in L \subset Y$ is called a \textit{focus} of a
generator $L$ if $dF \in L$ under some condition on the basis
forms $\theta^p$. In order to find the foci, we represent them in
the form $F = x^i A_i$. Then
$$
dF \equiv x^i \omega_i^p A_p \pmod{L},
$$
and as a result, the foci are determined by the system of
equations
$$
x^i \omega_i^p = 0.
$$
By (15), this system takes the form
\begin{equation}\label{16}
  x^i c_{iq}^p \theta^q = 0.
\end{equation}
This system has a nontrivial solution with respect to the forms
$\theta^q$ if and only \nolinebreak if
\begin{equation}\label{17}
  \det (x^i c_{iq}^p) = 0.
\end{equation}
Equation (17) determines on $L$ the focus hypersurface $F_L$ which
is an algebraic hypersurface of degree $r$.

Suppose that the point $A_0$ of our moving frame does not belong
to the hypersurface $F_L$. Then the 1-forms $\omega_0^p$ are
linearly independent, and we can take these forms as basis forms
of the congruence $Y$. As a result, equations (15) become
\begin{equation}\label{18}
  \omega_a^p = c_{aq}^p \omega_0^q,
\end{equation}
where $a = 1, \ldots , l,$ and $c_{0q}^p = \delta_q^p$. Now
equations (18) coincide with equations (6). As a result, equation
(17) of the focus hypersurface $F_L$ takes the form
\begin{equation}\label{19}
  \det (x^0 \delta_q^p + x^a c_{aq}^p) = 0.
\end{equation}
Equation (19) coincides with equation (12) defining the foci on a
plane generator $L$ of a variety $X$ with a degenerate Gauss map
of rank $r$. However, unlike in Section 1, the quantities
$c_{aq}^p$ are not connected by any relations of type (7), because
now there is no matrices $B^\alpha = (b_{pq}^\alpha)$. Thus, the
focus hypersurfaces $F_L$ determined by equation (19) are
arbitrary determinant varieties on generators $L$ of the
congruence $Y$ in question.

In particular, if $l = 1$ and $n = r + 1$, then $Y$ becomes a
rectilinear congruence. Equation (19) defining the focus
hypersurfaces $F_L$ of such a congruence becomes
\begin{equation}\label{20}
  \det (x^0 \delta_q^p + x^1 c_{1q}^p) = 0.
\end{equation}
Hence, each of the focus hypersurfaces $F_L$ of $Y$ decomposes
into $r$ real or complex points if each is counted as many times
as its multiplicity. Each of these points describes a
\textit{focal variety} in $\mathbb{P}^n$ tangent to the generators
$L$ of the congruence $Y$.

Next, we consider a pseudocongruence $Y^*$ in the space
$\mathbb{P}^n$. Its generator $L^*$ is of dimension $r - 1$ and
depends on $r$ parameters. We place the points $A_p, \linebreak p
= l + 1, \ldots , n, \, l = n - r$, of our moving frame into the
generator $L^* \subset Y^*$ and place the points $A_i, \, i = 0,
1, \ldots , l$, outside of $L^*$. The equations of infinitesimal
displacement of such frames again have the form (14) but now the
1-forms $\omega_p^i$ are linear combinations of the basis forms
$\theta^p$ defining a displacement of the generator $L^* = A_{l+1}
\wedge \ldots \wedge A_n$. So now we have
\begin{equation}\label{21}
  \omega_p^i = b_{pq}^i \theta^q
\end{equation}
and
\begin{equation}\label{22}
 dA_p = \omega_p^q A_q + b_{pq}^i \theta^q A_i.
\end{equation}

Consider a hyperplane $\xi$ passing through the generator $L^*
\subset Y^*$. Relative to our moving frame, the equation of $\xi$
is $\xi_i x^i = 0$, where $\xi_i$ are tangential coordinates of
the hyperplane $\xi$. The hyperplane $\xi$, which in addition to
the generator $L^*$ contains also a near generator ${}'L^*$
determined by the points $A_p$ and $dA_p$, is called the
\textit{focus hyperplane}. By (22), the conditions defining the
focus hyperplane are
\begin{equation}\label{23}
\xi_i b_{pq}^i \theta^q = 0.
\end{equation}
The system of equations (23) defines a displacement of the
generator $L^*$ if and only if this system has a nontrivial
solution with respect to the forms $\theta^q$. The necessary and
sufficient condition for existence of such a nontrivial solution
is the vanishing of the determinant of system (23):
\begin{equation}\label{24}
\det (\xi_i b_{pq}^i) = 0.
\end{equation}

Equation (24) defines the family of focus hyperplanes passing
through the generator $L^* \subset Y^*$. This family is an
algebraic hypercone of degree $r$ whose vertex is the generator
$L^*$. Note that equation (24) is similar to equation (13) of the
focus hypercone $\Phi_L$ of a variety with a degenerate Gauss map.

\section{Normalized Varieties in a Multidimensional
Projective Space} \textbf{1.} Consider a smooth $r$-dimensional
variety $X$ in a projective space $\mathbb{P}^n, \, r < n$. The
differential geometry on such a variety is rather poor. It is less
rich than the differential geometry on varieties of the Euclidean
space $\mathbb{E}^n$ or the spaces of constant curvature
$\mathbb{S}^n$ and $\mathbb{H}^n$, where by $\mathbb{S}^n$ and
$\mathbb{H}^n$ we denote, respectively, an $n$-dimensional
elliptic and hyperbolic space. With a first-order neighborhood of
a point $x \in X \subset \mathbb{P}^n$, only the tangent subspace
$T_x (X)$ is associated. For example, in Section 1.4 of the book
[AG 04], the authors showed that in order to enrich the
differential geometry of a curve in the projective plane
$\mathbb{P}^2$, it is necessary to use differential prolongations
of rather higher orders of the curve equations.

However, we can enrich the differential geometry of $X \subset
\mathbb{P}^n$ if we endow $X$ with an additional construction
consisting of a subspace $N_x (X)$ of dimension $n - r$ such that
$T_x (X) \cap N_x (X) = x$, and an $(r - 1)$-dimensional subspace
$K_x (X), \, K_x (X) \subset T_x (X), \, x \notin K_x (X)$. We
shall denote these subspaces simply by $N_x$ and $K_x$ and call
the \textit{normals of the first and second kind} (or simply the
\textit{first} and \textit{second normals}) of the variety $X$,
respectively (see [N 76], p. 198). The family of first normals
forms a \textit{congruence} $N$, and the family of second normals
forms a \textit{pseudocongruence} $K$ in the space $\mathbb{P}^n$.
If at any point of $x \in X$, there are assigned a single first
normal $N_x$ and a single second normal $K_x$, then the variety
$X$ is called \textit{normalizated} (cf. [N~76], p. 198, and [AG
93], Chapter 6).

As we will see below, on varieties of the Euclidean space
$\mathbb{E}^n$ and the non-Euclidean spaces $\mathbb{S}^n$ and
$\mathbb{H}^n$, the first and second normals are determined by the
geometry of these spaces while on varieties of the affine space
$\mathbb{A}^n$ and the projective space $\mathbb{P}^n$, these
normals should be assigned artificially, or to find them, one
should use higher order neighborhoods of a point $x \in X$. In
this section, we shall apply the first method. Note that the
second method is connected with great computational difficulties.
One can find more details on this method and a related
bibliography in
the books [AG 93], Chapters 6, 7, and [N 76], Chapter 5.

Thus, we now consider a normalized variety $X$ of dimension $r, \,
r = \dim X$, in the projective space $\mathbb{P}^n$. We associate
with $X$ a family of projective frames $\{A_0, A_1, \ldots ,
A_n\}$ in such a way that $A_0 = x; A_a \in N_x, \, a = 1, \ldots,
l$, where $l = n - r$; and $A_p \in K_x, \, p = l + 1,
 \ldots, n$. The equations of
infinitesimal displacement of these frames have the form
\begin{equation}\label{eq:25}
\renewcommand{\arraystretch}{1.5}
\left\{
\begin{array}{ll}
dA_0 = \omega_0^0 A_0 &\!\!\!\!+ \omega^p A_p , \\
dA_a = \omega_a^0 A_0 +  \omega_a^b A_b
       &\!\!\!\!+ \omega^p_a A_p, \\
dA_p = \omega_p^0 A_0 + \omega^a_p A_a &\!\!\!\!+ \omega_p^q A_q,\\
\end{array}
\right.
\renewcommand{\arraystretch}{1}
\end{equation}

Equations (25) show that for the family of moving frames in
question, the  system of differential equations
\begin{equation}\label{eq:26}
\omega^a = 0
\end{equation}
is satisfied, and the 1-forms $\omega^p$ are basis forms, because
they determine a displacement of the point $A_0 = x$ along the
variety $X$. Exterior differentiation of equations (26) and
application of Cartan's lemma lead to the following equations:
\begin{equation}\label{eq:27}
\omega_p^a = b_{pq}^a \omega^q, \;\; b_{pq}^a = b_{qp}^a
\end{equation}
The quantities $b^a_{pq}$ form a tensor and are coefficients of
the second fundamental forms of the variety $X$ at the point $x$
(see [AG 04], Section 2.1):
\begin{equation}\label{eq:28}
\Phi^a = b_{pq}^a \omega^p \omega^q.
\end{equation}

\textbf{2.} The points $A_p$ belong to the tangent subspace $T_x
(X)$. We assume that these points belong to the second normal $K_x
\subset T_x (X), \, K_x = A_{l+1} \wedge  \ldots \wedge A_n$.
Then, for $\omega^p = 0$, the 1-forms $\omega_p^0$ must also
vanish, and as a result, we have
\begin{equation}\label{eq:29}
\omega_p^0 = l_{pq} \omega^q.
\end{equation}

Next, we place the points $A_a$ of our moving frame into the first
normal $N_x$  of $X$, $N_x = A_0 \wedge A_{1} \wedge \ldots \wedge
A_l$. Then, for $\omega^p = 0$, we obtain that $\omega_a^p = 0$,
and hence
\begin{equation}\label{eq:30}
\omega_a^p = c_{a q}^p \omega^q.
\end{equation}

Consider a point $y \in N_x$ on the first normal. For this point,
we have \linebreak $y = y^0 A_0 + y^a A_a$. Differentiating this
point by means of (25), we find that
\begin{equation}\label{eq:31}
 dy = (dy^0 + y^0 \omega_0^0 +  y^a
\omega_a^0) A_0 + (y^0 \omega^p +   y^a \omega_a^p) A_p + (dy^a +
y^b  \omega_b^a) A_a.
\end{equation}

A point $y$ is a \textit{focus} of the first normal $N_x$ if $dy
\in N_x$. By (31), this condition implies that
$$
y^0 \omega^p +   y^a \omega_a^p = 0.
$$
Applying relations (30), we find that
$$
(y^0 \delta_q^p +  y^a c_{a q}^p) \omega^q = 0.
$$
This system has a nontrivial solution with respect to the forms
$\omega^q$ if and only \nolinebreak if
\begin{equation}\label{eq:32}
\det (y^0 \delta_q^p +  y^a c_{a q}^p) = 0.
\end{equation}
Equation (32) differs from equation (19) only in notation, and it
defines the focus hypersurface $F_x$ in the generator $N_x$ of the
congruence of first normals associated with the variety $X$. It
follows from equation (32) that the point $x \in X$, whose
coordinates are $y^0 = 1, \, y^a = 0$, does not belong to the
focus hypersurface $F_x$.

Let us find the focus hypercones $\Phi_x$ of the pseudocongruence
$K$ of second normals of $X$. The hypercones $\Phi_x$ are formed
by the hyperplanes $\xi$ of the space $\mathbb{P}^n$ containing
the second normal $K_x = A_{l+1} \wedge \ldots \wedge A_n \subset
T_x (X)$ and its neighboring normal $K_x + d K_x$, which contains
not only the points $A_p$ but also the points
$$
d A_p \equiv \omega_p^0 A_0 + \omega_p^a A_a \pmod{N_x}.
$$
As a result, tangential coordinates $\xi_0$ and $\xi_a$ of such a
hyperplane satisfy the equations
$$
\xi_0 \omega_p^0 +  \xi_a\omega_p^a = 0.
$$
By (29) and (27), it follows from this equation that
$$
(\xi_0 l_{pq} +  \xi_a b^a_{pq})\, \omega^q = 0.
$$
This system has a nontrivial solution with respect to the forms
$\omega^q$ if and only if its determinant vanishes,
\begin{equation}\label{eq:33}
\det (\xi_0 l_{pq} +  \xi_a b^a_{pq}) = 0.
\end{equation}
Equation (33) determines an algebraic hypercone of order $r$ whose
vertex is the generator $K_x$ of the pseudocongruence $K$ of the
second normals. This hypercone is called the \textit{focus
hypercone} of the pseudocongruence $K$.

\textbf{3.} Next, we consider the tangent and normal bundles
associated with a normalized variety $X$. The base of both bundles
is the variety $X$ itself, the fibers of the tangent bundle are
the tangent subspaces $T_x$, and  the fibers of the normal bundle
are the first normals $N_x$,

Suppose that ${}'x = x + x^p A_p$ is an arbitrary point in the
tangent subspace $T_x$, and $\boldsymbol{x} = {}'x - x =
 x^p A_p$ is a vector in the tangent bundle $TX$. The differential of
 this vector has the form
\begin{equation}\label{eq:34}
 d\boldsymbol{x} = (dx^p + x^q \omega_q^p) A_p + x^p (l_{pq} A_0 +  b^a_{pq} A_a) \omega^q.
\end{equation}
The first term on the right-hand side of (34) belongs to the
tangent subspace $T_x$, and the second term belongs to $N_x$. The
1-form
$$
D x^p = dx^p + x^q \omega_q^p
$$
is called the \textit{covariant differential} of the vector
$\boldsymbol{x} = (x^p)$ in the connection $\gamma_t$, and the
1-forms $\omega_q^p$ are the components of the \textit{connection
form} $\omega =\{\omega_q^p\}$ of the affine connection $\gamma_t$
induced on the variety $X$ by a normalization of $(N, K)$.

The vector field $\boldsymbol{x}$ is called \textit{parallel} in
the connection $\gamma_t$ if its covariant differential  $DX^p$
vanishes, i.e, if
\begin{equation}\label{eq:35}
Dx^p = dx^p + x^q \omega_q^p = 0.
\end{equation}

We find the exterior differentials of the components $\omega_q^p$
of the connection form $\omega$. By (27), (29), and (30), these
exterior differentials have the form
\begin{equation}\label{eq:36}
d \omega_q^p = \omega_q^s \wedge \omega_s^p + (l_{qs} \delta_t^p +
b^a_{qs} c_{at}^p) \omega^s \wedge \omega^t.
\end{equation}
The 2-form
$$
\Omega_q^p = d \omega_q^p - \omega_q^s \wedge \omega_s^p
$$
is said to be \textit{the curvature form} of the affine connection
$\gamma_t$ induced on  the variety $X$. From equation (36) it
follows that
\begin{equation}\label{eq:37}
\Omega_q^p = \frac{1}{2} R^p_{qst}  \omega^s \wedge \omega^t,
\end{equation}
where
\begin{equation}\label{eq:38}
R^p_{qst} = l_{qs} \delta_t^p + b^a_{qs} c_{at}^p - l_{qt}
\delta_s^p - b^a_{qt} c_{as}^p
\end{equation}
 is the \textit{curvature tensor} of the affine
connection $\gamma_t$ on $X$. Equations (38) allow us to compute
the curvature tensor for different normalizations of the variety
$X$.

If $R^p_{qst} = 0$ on the variety $X$, then the affine connection
$\gamma_t$ on $X$ is \textit{flat}, and a parallel translation of
a vector $\boldsymbol{x} $ does not depend on the path of
integration (see, for example, [N 76], p. 118, or [KN 76], p. 70).

\textbf{4.} Further, we consider a vector field $\boldsymbol{y}$
in the normal bundle $N (X)$. This vector is determined by the
point $x$ and a point $y = y^0 A_0 + y^a A_a$ of the fiber $N_x
\subset N (X)$. The differential of the point $y$ is defined by
equation (31).

The 1-form
\begin{equation}\label{eq:39}
D y^a = dy^a + y^b \omega_b^a
\end{equation}
is called the \textit{covariant differential} of the vector field
$\boldsymbol{y}$ in the normal bundle $N (X)$, and the forms
$\omega_a^b$ are the components of the \textit{connection form of
the normal connection} $\gamma_n$ on a normalized variety $X$
(see, for example, [Ca 01], p.~242; see more on the normal
connection in [AG 95] and Section 6.3 of the book [AG 93]). The
2-form
$$
\Omega_b^a = d \omega_b^a - \omega_b^c \wedge \omega_c^a
$$
is called the \textit{curvature form} of the normal connection
$\gamma_n$. Note that Cartan in [Ca~01] called this form the
\textit{Gaussian torsion} of an embedded variety $X$.

 Differentiating the forms $\omega_b^a$ and applying formulas
(27) and (30), we find the expression of the curvature form
$\Omega_b^a$:
\begin{equation}\label{eq:40}
\Omega_b^a = \frac{1}{2} R^a_{bst} \omega^s \wedge \omega^t,
\end{equation}
where
\begin{equation}\label{eq:41}
R^a_{bst} =  c_{bs}^p  b^a_{pt}-  c_{bt}^p b^a_{ps}.
\end{equation}
The tensor $R^a_{bst}$ is called the \textit{tensor of normal
curvature} of the variety $X$.

The second normals $K_x$ associated with the variety $X$ allow us
to find a distribution $\Delta_y$ of $r$-dimensional subspaces
associated with $X$. The elements of the distribution $\Delta_y$
are linear spans of the points $y \in N_x$ and the second normals
$K_x, \; \Delta_y = y \wedge K_x$. By (31),  the distribution
$\Delta_y$ is determined by the system of equations
\begin{equation}\label{eq:42}
dy^a + y^b \omega^a_b = 0.
\end{equation}
In the general case, the system of equations (42) is not
completely integrable, and when a point $x$ moves along a closed
contour $l \subset X$, the corresponding point $y$ does not
describe a closed contour.

But the point $y$ describes a closed contour $l'$ if system (42)
is completely integrable. The condition of complete integrability
of (42) is the vanishing of the tensor of normal curvature (41) of
the variety $X$. In this case, the distribution $\Delta_y$ defined
by system (42)   is completely integrable, and the closed contours
$l'$ lie on integral varieties of this distribution. These
integral varieties form an $(n - r)$-parameter family of
$r$-dimensional subvarieties $X (y)$ which are ``parallel'' to the
variety $X$ in the sense that the subspaces $T_x (X)$ and $T_x (X
(y))$ pass through the same second normal $K_x$.


\textbf{5.} A normalization of a variety $X \subset \mathbb{P}^n$
is called \textit{central} if all its first normals $N_x$ form a
bundle with an $(l-1)$-dimensional vertex $S$.

The following theorem gives necessary and sufficient conditions
for a normalization of a variety $X$ to be central.

\begin{theorem}
A normalization of a  normalized variety $X \subset \mathbb{P}^n$
is central if and only if the quantities $l_{pq}$ and $c^p_{aq}$
in equations $(29)$ and $(30)$ vanish:
$$
l_{pq} = 0, \;\; c^p_{aq} = 0.
$$
\end{theorem}

\begin{proof} {\sf Necessity}:
 Suppose that a normalization of a variety $X \subset
\mathbb{P}^n$ is central with an $(l-1)$-dimensional vertex $S$.
If we place the points $A_a$ into this vertex $S$, then we get
$$
d A_a = \omega_a^b A_b.
$$
By (25), this implies that
$$
\omega_a^0 = 0, \;\; \omega_a^p = 0.
$$
By (29) and (30), this means that
$$
l_{pq} = 0, \;\; c^p_{aq} = 0.
$$

{\sf Sufficiency}: If $l_{pq} = 0, \; c^p_{aq} = 0,$ then $
\omega_a^0 = 0, \; \omega_a^p = 0,$ and $d A_a = \omega_a^b A_b, $
(i.e., the subspace $S = A_{1} \wedge \dots \wedge A_l$ is fixed),
then all first $l$-dimensional first normals $N_x$ pass through
$S$, and a normalization of $X$ is central with the
$(l-1)$-dimensional vertex $S$.
\end{proof}

\begin{corollary}
The induced affine connection $\gamma_t$ and the normal connection
$\gamma_n$ of a centrally normalized variety $X \subset
\mathbb{P}^n$ are flat.
\end{corollary}

\begin{proof} Because for a centrally normalized  variety $X$,
we have $l_{pq} = 0, \; c^p_{aq} = 0,$ and the curvature tensor of
the induced affine connection $\gamma_t$ has the form (38), it
follows that
$$
R^p_{qst} = 0,
$$
i.e.,  the curvature tensor of the induced affine connection
$\gamma_t$ of a centrally normalized variety vanishes.

In the same way, it follows from (41) that  the  tensor of the
normal connection $\gamma_n$ of a centrally normalized variety
also vanishes.

Note that both results also follow from the fact that a centrally
normalized variety $X \subset \mathbb{P}^n$ can be bijectively
projected onto an $r$-dimensional subspace $T$ that is
complementary to the vertex $S$ of the bundle of first normals
$N_x$, and the geometry on the variety $X$ induced by this central
normalization is equivalent to the plane geometry in the subspace
$T$.
\end{proof}

Atanasyan [A 52] found necessary and sufficient conditions for
normalization of a variety $X$ in an affine space $\mathbb{A}^{N}$
to be central and trivial. The trivial normalization  of $X$ in
$\mathbb{A}^N$ is a normalization for which all first normals
$N_x$ of $X$ are parallel to some constant $l$-plane (i.e., they
form a bundle of parallel $l$-planes).

In our notations, his conditions for a normalization to be central
are
$$
c_{aq}^p = \delta_q^p c_a,
$$
where $\delta_q^p$ is the Kronecker delta, and $c_a$ are $(1,
0)$-tensors, and the conditions for a normalization to be trivial
are
$$
c_{aq}^p = 0.
$$
But for an affine space (and in particular,
for a Euclidean space) we always have $l_{pq} = 0$. In addition,
in the projective setting (as well as in the  affine setting), a
trivial normalization is a central normalization whose vertex is
an $(l-1)$-plane at infinity. Therefore, Atanasyan's results
follow from Theorem 1.

\textbf{6.} Consider the normalization dual to a central
normalization. For such a normalization all second normals $K_x$
belong to a fixed hyperplane $\alpha$. We will call such a
normalization \textit{affine}.

\begin{theorem} A normalization of a variety $X \subset
\mathbb{P}^n$ is affine if and only if the $1$-forms $\omega_p^0$
and $\omega_a^0$ occurring in equations $(25)$ vanish,
$$
\omega_p^0 = 0, \;\; \omega_a^0 = 0.
$$
If a normalization of variety $X \subset \mathbb{P}^n$ is affine,
then the space $\mathbb{P}^n$ carries an affine structure, i.e.,
$\mathbb{P}^n$ is an affine space $\mathbb{A}^n$
\end{theorem}

\begin{proof}
We place the points $A_1, \dots, A_n$ of our moving frame into the
fixed hyperplane $\alpha$. Since for an affine normalization $K_x
\subset \alpha$, and hence $dA_p \subset \alpha, \linebreak p =
l+1, \dots, n$, it follows that
$$
\omega_p^0 = 0.
$$

Moreover, the points $A_a, \, a = 1, \dots, l$, of the first
normal $N_x$ can be also placed into the hyperplane $\alpha$. Then
$dA_a \subset \alpha$, and as a result, we have
$$
\omega_a^0 = 0.
$$

Conversely, if $\omega_p^0 = 0, \omega_a^0 = 0$, then
$$
dA_p \subset \alpha, \;\; d A_a \subset \alpha,
$$
where $\alpha = A_1 \wedge \dots \wedge A_n$. Hence the plane
$\alpha$ is fixed, and the normalization of $X$ is affine. This
proves the first part of Theorem 3.

To prove the second part of Theorem 3, note that we can take the
hyperplane $\alpha$ as the hyperplane at infinity
$\mathbb{P}^{n-1}_\infty$ of the space $\mathbb{P}^n$. Hence, this
hyperplane defines an affine structure of the space
$\mathbb{P}^n$. Thus, the space $\mathbb{P}^n$ become an affine
space $\mathbb{A}^n$.
\end{proof}

\textbf{7.} Now suppose that a normalized variety $X \subset
\mathbb{P}^n$ has a flat normal connection $\gamma_n$, i.e.,
$R^a_{bst} = 0$. By (41), these conditions lead to the relation
\begin{equation}\label{eq:43}
b^a_{pt} c^p_{bs} = b^a_{ps} c^p_{bt}.
\end{equation}
Relations (43) differ from relations (7) only in notation. If we
introduce the matrix notations
$$
B^a = (b^a_{pq}), \;\;  C_b = (c^p_{bq}),
$$
 then relations (43) take the form
\begin{equation}\label{eq:44}
(B^a C_b) = (B^a C_b)^T.
\end{equation}

We proved in Chapters 3 and 4 of [AG 04] that these relations
imply that the matrices $B^a$ and $C_b$ can be simultaneously
reduced to a diagonal form or a block diagonal form. Thus, we have
proved the following result.

\begin{theorem}
 The focus hypersurfaces $F_x \subset N_x$ of a normalized variety $X$
with a flat normal connection decompose into the plane generators
of different dimensions.
\end{theorem}

This property of the varieties $X$ with a flat normal connection
$\gamma_n$ allows us to construct a classification of such
varieties in the same way as this was done for the varieties with
degenerate Gauss maps in a projective space. For varieties in an
affine space and a Euclidean space, such a classification was
outlined in the papers [ACh 75, 76, 01].

\section{Normalization of Varieties in Affine \\ and Euclidean Spaces}

\textbf{1.} An affine space $\mathbb{A}^n$ differs from a
projective space $\mathbb{P}^n$ by the fact that in $\mathbb{A}^n$
a hyperplane at infinity $\mathbb{P}_\infty$ is fixed. If we place
the points $A_i, \, i = 1, \dots, n$, of our moving projective
frame into this hyperplane, then the equations of infinitesimal
displacement of the moving frame take the following form (see
equations (1.81) in {AG 04]):
\begin{equation}\label{eq:45}
\renewcommand{\arraystretch}{1.3}
\left\{
\begin{array}{lll}
dA_0 = & \!\!\!\! \omega^0_0 A_0 + & \!\!\!\!  \omega^i_0 A_i, \\
dA_i = & \!\!\!\!  & \!\!\!\!  \omega^j_i A_j, \; i, j = 1, \ldots
, n,
\end{array}
\right.
\renewcommand{\arraystretch}{1}
\end{equation}
 and the structure equations of the  affine space $\mathbb{A}^n$
take the form
\begin{equation}\label{eq:46}
d \omega^0_0 = 0, \;\; d \omega^i_0 = \omega_0^j \wedge
\omega_j^i, \;\; d \omega^i_j = \omega_j^k \wedge \omega_k^i.
\end{equation}

Consider a variety $X$ of dimension $r$ in the affine space
$\mathbb{A}^n$. The tangent space $T_x (X)$ intersects the
hyperplane at infinity $\mathbb{P}_\infty$ in a subspace $K_x$ of
dimension $r - 1$, $K_x = T_x \cap \mathbb{P}_\infty$. Thus, for a
normalization of $X$, it is sufficient to assign only a family of
first normals $N_x$. If we place the points $A_a, \linebreak a =
1, \dots , l$, of our moving frame into the subspace $N_x \cap
\mathbb{P}_\infty$, and the points $A_p, \,p = l+1, \dots , n$,
into the subspace $K_x$, then equations (45) take the form
\begin{equation}\label{eq:47}
\renewcommand{\arraystretch}{1.3}
\left\{
\begin{array}{llll}
dA_0 = &\!\!\!\! \omega^0_0 A_0 + & & \!\!\!\! + \omega^p_0 A_p, \\
dA_a = & &\!\!\!\!  \omega_a^b A_b  & \!\!\!\! +\omega^p_a A_p, \\
dA_p = & &\!\!\!\!  \omega_p^a A_a  &\!\!\!\! + \omega^q_p A_q
\end{array}
\right.
\renewcommand{\arraystretch}{1}
\end{equation}
(cf. equations (25)).

As was in the projective space, we have equations (27),
\begin{equation}\label{eq:48}
\omega_p^a = b_{pq}^a \omega^q, \;\; b_{pq}^a = b_{qp}^a,
\end{equation}
where $ b_{pq}^a$ is the second fundamental tensor of the variety
$X$. Equations (30) also preserve their form:
\begin{equation}\label{eq:49}
\omega_a^p = c_{a q}^p \omega^q,
\end{equation}
but equations (29) become
\begin{equation}\label{eq:50}
\omega_p^0 = 0.
\end{equation}

As a result, because $l_{pq} = 0$, the equation of the focus
hypersurface $F_x \subset N_x$ preserves its form (32):
\begin{equation}\label{eq:51}
\det (y^0 \delta_q^p +  y^a c_{a q}^p) = 0.
\end{equation}
As to  equation (33) of the focus hypercone $\Phi_x$, by (50),
this equation takes the form
\begin{equation}\label{eq:52}
\det (\xi_a b^a_{pq}) = 0.
\end{equation}

Expressions (38) for the components of the curvature tensor of the
affine connection $\gamma_t$  induced on the normalized variety $X
\subset \mathbb{A}^n$ take now the form
\begin{equation}\label{eq:53}
R^p_{qst} = b^a_{qs} c_{at}^p - b^a_{qt} c_{as}^p,
\end{equation}
and the expression (41) for the components of the tensor of normal
curvature of the variety $X$ preserves its form:
\begin{equation}\label{eq:54}
R^a_{bst} =  b^a_{pt} c_{bs}^p -  b^a_{ps} c_{bt}^p.
\end{equation}

Consider the tensor $R_{st}$ obtained from the curvature tensor
$R^p_{qst}$ of the affine connection $\gamma_t$ with respect to
the indices $p$ and $q$. This tensor is called the
\textit{Ricci-type tensor of the connection $\gamma_t$}. It
follows from (53) that
$$
R_{st} = b^a_{ps} c^p_{at} - b^a_{pt} c^p_{as}.
$$

Similarly we can define the \textit{Ricci-type tensor of the
normal connection} $\gamma_n$. We denote it by
$\widetilde{R}_{st}$. It follows from (54) that
$$
\widetilde{R}_{st} = b^a_{pt} c^p_{as} - b^a_{ps} c^p_{at}.
$$
Comparing the last two equations, we see that
$$
R_{st} = -\widetilde{R}_{st}.
$$

Hence, the following theorem is valid.

\begin{theorem}
On a variety $X \subset \mathbb{A}^n$ endowed with an affine
normalization, the Ricci-type tensors of the connections
$\gamma_t$ and $\gamma_n$ are equal in absolute value but opposite
in sign.
\end{theorem}

Note that if a normalized variety $X$ is a hypersurface in the
space $\mathbb{A}^n$, then the following theorem is valid.

\begin{theorem} If on a normalized hypersurface $X \subset \mathbb{A}^n$
the induced affine connection
$\gamma_t$ is flat, then the normal connection $\gamma_n$ is also
flat.
\end{theorem}

\begin{proof} In fact, for a hypersurface $X$ we have the
following ranges of the indices:
$$
a, b = 1; \;\; p, q, s, t = 2, \dots , n.
$$
Hence the curvature tensor of the normal connection $\gamma_n$ has
the components $R^1_{1st}$. Now it follows from Theorem 5 that
$$
R^1_{1st} = - R^p_{pst}.
$$
But if the  connection $\gamma_t$ is flat, then $R^p_{qst} = 0$
and hence $R^p_{pst} = 0$. As a result, we have $R^1_{1st} = 0$,
and the  connection $\gamma_n$ is also flat.
\end{proof}

As was the case in the projective space, the vanishing of tensor
of normal curvature $R^a_{bst}$ is equivalent to the complete
integrability of the system defining the distribution $\Delta_y =
y \wedge K_x$, where $y \in N_x$. But in the affine space, the
elements $\Delta_y$ of this distribution are parallel to the
subspace $T_x (X)$.

Thus, we have proved the following result.

\begin{theorem}
A variety $X \subset \mathbb{A}^n$ has a flat normal connection
$\gamma_n$ if and only if this variety admits an $l$-parameter
family of parallel varieties $X (y)$, where $y \in N_x$.
\end{theorem}

\textbf{2.}  Another relation of the theory of varieties with
degenerate Gauss maps and the theory of normalized varieties was
established by Chakmazyan in [Cha 77]  (see also Theorem 4 in [Cha
78] and Theorem 1 on p. 39 of his book [Cha 90]).

We will now present this theorem.

Suppose that at an arbitrary point $x$ of a normalized variety
$X$, an $s$-dimensional direction $\nu^s (x)$ (i.e., an a
$s$-dimensional plane passing through $x$) belonging to the first
normal $N_x (X)$ is given. This means we have a smooth field of
$s$-dimensional directions $\nu^s (x)$ on $X$, where $s \leq l = n
- r$. This field determines a normal subbundle $\nu^s (X)$, whose
$s$-dimensional fibres are $s$-dimensional centroprojective
spaces.

The plane $N^s (x)$ of this field corresponding to an arbitrary
point $x \in X$, can be defined by the point $x$ and points $B_f$
given by
\begin{equation}\label{eq:55}
B_f = \xi_f^a A_a \in N_x,
\end{equation}
where $f, g, h = 1, \dots, s$.

In addition, the plane $N^s (x)$ must be invariant with respect to
admissible transformations of the moving frame in $N_x (X)$. The
necessary and sufficient conditions for this invariance are
\begin{equation}\label{eq:56}
dB_f = \theta_f^g B_g + \theta_f^0 A_0 \pmod{\omega^p},
\end{equation}
where $\theta_f^g $ and $\theta_f^0$ are linearly independent
1-forms whose structure can be obtained by taking exterior
derivatives of (56). Since we do not need these conditions, we
will not derive them.

A field $\nu^s$ is called \textit{parallel} with respect to the
normal connection $\gamma_n$ if under any infinitesimal
displacement of an arbitrary point $x \in X$, the displacement of
the $s$-dimensional direction $\nu^s (x)$ remains in the $(r +
s)$-dimensional plane passing through the tangent subspace $T_x
(X), \, x \in X$, and the direction $\nu^s (x)$.

Let us find analytic conditions for a field $\nu^s$ to be
parallel. Any direction belonging to an $s$-dimensional element of
the field $\nu^s$ is determined by the point $A_0 = x$ and
\begin{equation}\label{eq:57}
A =  A_0 + \xi^f B_f,
\end{equation}
where $B_f$ are defined by (55).

Taking exterior derivative of (57) and applying (25), we find that
\begin{equation}\label{eq:58}
dA = (\omega_0^0 + \xi^f \xi_f^a \omega_a^0) A_0 + d \xi^f B_f +
(\omega^p + \xi^f \xi_f^a \omega_a^p) A_p + \xi^f (d\xi_f^a +
\xi_f^b \omega_b^a) A_a.
\end{equation}
The field $\nu^s$ is parallel with respect to the normal
connection $\gamma_n$ if and only if
$$
(d\xi_f^a + \xi_f^b \omega_b^a) A_a = \theta_f^g B_g + \theta_f^0
A_0.
$$
This and (55) implies that
\begin{equation}\label{eq:59}
 d\xi_f^a + \xi_f^b \omega_b^a = \theta^g_f \xi_g^a.
\end{equation}

\begin{theorem} A field $\nu^s$ of $s$-dimensional
normal directions $\nu^s (x)$ on a normalized variety $X \subset
\mathbb{P}^n$ is parallel with respect to the normal connection
$\gamma_n$ if and only if all the planes $N^s (x),\, x \in X$,
form a variety $V^{r+s}_{r}$ with a degenerate Gauss map of rank
$r$ with $s$-dimensional plane generators in the space
$\mathbb{P}^n$.
\end{theorem}

\begin{proof} See [Cha 77] or [Cha 90], pp. 39-40. \end{proof}



Theorem 8 indicates a method of construction of a  variety
$V^{r+s}_{r}$ with a degenerate Gauss map departing from a
normalized variety $X$ of some special kind.

\textbf{3.} Further consider a variety $X$ of dimension $r$ in the
Euclidean space $\mathbb{E}^n$. On $X$, both the second normal
$K_x = T_x \cap \mathbb{P}_\infty$ and the first normal $N_x$
orthogonal to the tangent subspace $T_x (X)$ are naturally
defined.

In the Euclidean space $\mathbb{E}^n$, there is defined a scalar
product of vectors, and a scalar product of points in the
hyperplane at infinity $\mathbb{P}_\infty$ is induced by the
scalar product in $\mathbb{E}^n$. Because in our moving frame, we
have  $A_a \in N_x \cap \mathbb{P}_\infty; \linebreak A_p \in T_x
\cap \mathbb{P}_\infty = K_x, \, a = 1, \dots, l; \, p = l + 1,
\dots, n$; and $T_x \perp N_x$, we find that
\begin{equation}\label{eq:60}
(A_a, A_p) = 0,
\end{equation}
where, as usually, the parentheses denote the scalar product of
points in the hyperplane at infinity $\mathbb{P}_\infty$. In
addition, we set
\begin{equation}\label{eq:61}
(A_a, A_b) = g_{ab}, \;\; (A_p, A_q) = g_{pq},
\end{equation}
where $g_{ab}$ and $g_{pq}$ are nondegenerate symmetric tensors.

Differentiating equations (60) and using formulas (47), (60), and
(61), we find that
$$
g_{ab}\, \omega_p^b + g_{pq}\, \omega_a^q = 0.
$$
It follows that
\begin{equation}\label{eq:62}
\omega_a^p = -g^{pq}\, g_{ab}\, \omega_q^b.
\end{equation}
Equations (57) and (48) imply that
\begin{equation}\label{eq:63}
\omega_a^p = -g^{pq}\, g_{ac} \, b_{qs}^c \omega^s.
\end{equation}
Comparing (63) and (49), we obtain
\begin{equation}\label{eq:64}
c_{as}^p = -g^{pq} \,g_{ac} \,b_{qs}^c.
\end{equation}

Now we find the equation of the focus hypersurface $F_x$ of the
variety \linebreak $X \in \mathbb{E}^n$. By (51) and (64), we have
the following equation for $F_x$:
$$
\det (y^0 \delta_q^p -  y^a g^{ps} g_{ac} b_{sq}^c) = 0.
$$
The last equation is equivalent to the equation
\begin{equation}\label{eq:65}
\det (y^0 g_{pq} -  y_a b_{p q}^a) = 0,
\end{equation}
where $y_a = g_{ab} y^b$.

In our moving frame, the hyperplane at infinity
$\mathbb{P}_\infty$ is determined by the equation $y^0 = 0$. Hence
by (65), the intersection $F_x \cap \mathbb{P}_\infty$ of the
focus hypersurface $F_x$ with the hyperplane at infinity
$\mathbb{P}_\infty$ is defined by the equation
\begin{equation}\label{eq:66}
\det (y_a b_{p q}^a) = 0.
\end{equation}
But this equation differs only in notation from equation (52) of
the focus hypercone $\Phi_x$ of the variety $X$. Equations (52)
and (66) coincide if $\xi_a = y_a = g_{ab} y^b$. Thus, we have
proved the following result.

\begin{theorem}
The focus hypercone $\Phi_x$ of the variety $X \subset
\mathbb{E}^n$ is formed by the hyperplanes $\xi$ containing the
tangent subspace $T_x$ and orthogonal at the points
$\widetilde{y}$ of the hyperplane at infinity $\mathbb{P}_\infty$
lying in the intersection $F_x \cap \mathbb{P}_\infty$.
\end{theorem}

This result clarifies the geometric meaning of the focus hypercone
$\Phi_x$ for the variety $X \subset \mathbb{E}^n$ and its relation
with the focus hypersurface $F_x$ of $X$.

We also find the curvature tensor  of the affine connection
induced on the variety $X \subset \mathbb{E}^n$. Substituting the
values of $c_{aq}^p$ from (64) into formula (53), we find that
\begin{equation}\label{eq:67}
R^p_{qst} = g^{pu} g_{ac} (b^a_{qt} b^c_{us} - b^a_{qs} b_{ut}^c).
\end{equation}
Contracting equation (67) with the tensor $g_{pv}$ and changing
the summation indices (if necessary), we find that
\begin{equation}\label{eq:68}
R_{pqst} =  g_{ac} (b^a_{ps} b^c_{qt} - b^a_{pt} b_{qs}^c),
\end{equation}
where $R_{pqst} = g_{pu} R^u_{qst}$. Formulas (67) and (68) give
the usual expressions for the curvature tensor of the affine
connection $\gamma_t$ induced on a normalized variety $X \subset
\mathbb{E}^n$.

But in addition to the curvature tensor  of the affine connection
induced on a normalized variety $X \subset \mathbb{E}^n$, we
considered also the tensor $R^a_{bst}$ of normal curvature defined
by equation (54). Substituting the values of $c_{aq}^p$ from (64)
into formula (54), we find that
\begin{equation}\label{eq:69}
R^a_{bst} = g^{pq} g_{bc} (b^c_{qt} b_{ps}^a -  b^c_{qs}
b_{pt}^a).
\end{equation}
As we noted earlier, in the book [C 01], the exterior 2-form
$$
\Omega_b^a = d \omega_b^a - \omega_b^c \wedge \omega_c^a =
\frac{1}{2} R^a_{bst} \omega^s \wedge \omega^t
$$
is called the \textit{Gaussian torsion} of a variety $X \subset
\mathbb{E}^n$.


\vspace*{10mm}

\noindent {\em Authors' addresses}:\\

\noindent
\begin{tabular}{ll}
M.~A. Akivis &V.~V. Goldberg\\
Department of Mathematics &Department of Mathematical Sciences\\
Jerusalem Institute of Technology---Mahon Lev &  New
Jersey Institute of Technology \\
Havaad Haleumi St., P. O. B. 16031 & University Heights \\
Jerusalem 91160, Israel &  Newark, N.J. 07102, U.S.A. \\

 E-mail address: akivis@mail.jct.ac.il & E-mail address:
 vlgold@m.njit.edu
 \end{tabular}

\vspace*{3mm}

 \noindent
A.~V.~Chakmazyan\\
Department of Geometry\\
Armyanskij Gosudarstvennyj Pedagogicheskij Institut \\
Erevan, Armenia

\end{document}